\documentclass[12pt]{article}

\usepackage{sectsty}
\sectionfont{\color{dark}} 
\subsectionfont{\color{dark}} 
\subsubsectionfont{\color{light}} 
\usepackage{fullpage}

\usepackage{cancel,enumerate}
\usepackage[affil-it]{authblk}

\usepackage[intlimits]{amsmath}
\usepackage{amstext, amssymb}
\usepackage{eufrak}
\usepackage[mathscr]{euscript}

\usepackage{tikz}[2013/12/13]
\usetikzlibrary{calc}
\usetikzlibrary{matrix}
\definecolor{dark}{RGB}{0, 133, 202 }
\definecolor{light}{RGB}{112, 155, 230}

\usepackage {
	amsmath,
	amssymb,
	amsthm,
	amsxtra,
	appendix,
	epsfig,
	fourier-orns,
	graphics,
	longtable,
	needspace,
	makeidx,
	mdwlist,
	multicol,
	psfrag,
	pstricks,
	rotating,
	setspace,
	sidecap,
	stmaryrd,
	subfloat,
	tabularx,
	textcomp,
	tocloft,
	url,
	wrapfig
}

\usepackage{etex}
\usepackage{mathtools}
\usepackage{needspace}

\allowdisplaybreaks

\AtBeginDocument{%
    \mathchardef\mathcomma\mathcode`\,
    \mathcode`\,="8000 
    \mathchardef\mathsemicolon\mathcode`\;
    \mathcode`\;="8000 
}
{\catcode`,=\active
  \gdef,{\mathcomma\discretionary{}{}{}}
\catcode`;=\active
  \gdef;{\mathsemicolon\discretionary{}{}{}}
}

\theoremstyle{definition}

\newtheorem{theorem}{Theorem}[section]

\newcommand{\pdiff}[2]{\frac{\partial #1}{\partial #2}}

\makeatletter
\newcommand{\superimpose}[2]{%
  {\ooalign{$#1\@firstoftwo#2$\cr\hfil$#1\@secondoftwo#2$\hfil\cr}}}
\makeatother

\renewcommand{\AA}{\mathbb{A}}
\renewcommand{\det}[1]{\left| #1 \right|}

\begin{document}

\title{\color{dark} Yet Another Proof of Sylvester's Determinant Identity}
\author[1]{Jan Vrbik}
\author[2]{Paul Vrbik}
\affil[1]{Department of Mathematics and Statistics\\ Brock University, Canada}
\affil[2]{School of Mathematical and Physical Sciences, University of Newcastle, Australia}
\date{\today}
\maketitle

\begin{abstract}

In 1857 Sylvester stated a result on determinants without proof that was recognized as important over the subsequent century. Thus it was a surprise to Akritas, Akritas and Malaschonok when they found only one English proof --- given by Bareiss 111 years later! To rectify the gap in the literature these authors collected and translated six additional proofs: four from German and two from Russian \cite{akritas}. These proofs range from long and ``readily understood by high school students''  to  elegant but high level.
\medskip

We add our own proof to this collection which exploits the product rule and the fact that taking a derivative of a determinant with respect to one of its elements yields its cofactor. A differential operator can then be used to replace one row with another.

\end{abstract}

\section{Preliminaries}
Let $\AA$ be an $n \times n$ matrix with entries $a_{i,j}$ for $1 \leq i,j \leq n$ and denote by $\AA_{i|k}$ the matrix $\AA$ with Row $i$ and Column $k$ deleted. Similarly let $\AA_{i,j|k,\ell }$ be the matrix obtained by deleting Rows $i$ and $j$ \emph{and} Columns $k$ and $\ell$:
\begin{center}
\hfill
\raisebox{6.8ex}{$\AA_{i|k}:=$}
\begin{tikzpicture}
    \matrix (M)[matrix of math nodes,left delimiter={[},right delimiter={]}]{
      ~ & a_{k,1} & ~ & ~ & ~ & ~\\
      ~ & ~ & ~ & ~ & ~ & ~\\
      a_{i,1} & ~ & ~ & ~ & ~ & a_{i,n} \\
      ~ & ~ & ~ & ~ & ~ & ~\\
      ~ & ~ & ~ & ~ & ~ & ~\\
      ~ & a_{k,n} & ~ & ~ & ~ & ~\\
     };
     \draw[black](M-3-1.west)--(M-3-6.east);
     \draw[black](M-1-2.north)--(M-6-2.south);
\end{tikzpicture}
\hspace{-.3em}\raisebox{6.8ex}{,}\;\;
\hfill
\raisebox{6.8ex}{$\AA_{i,j|k,\ell}:=$}
\begin{tikzpicture}
    \matrix (M)[matrix of math nodes,left delimiter={[},right delimiter={]}]{
      ~ & a_{k,1} & ~ & a_{\ell,1} & ~ & ~\\
      ~ & ~ & ~ & ~ & ~ & ~\\
      a_{i,1} & ~ & ~ & ~ & ~ & a_{i,n} \\
      ~ & ~ & ~ & ~ & ~ & ~\\
      a_{j,1} & ~ & ~ & ~ & ~ & a_{j,n} \\
      ~ & a_{k,n} & ~ & a_{\ell,n} & ~ & ~\\
     };
     \draw[black](M-3-1.west)--(M-3-6.east);
     \draw[black](M-5-1.west)--(M-5-6.east);
     \draw[black](M-1-2.north)--(M-6-2.south);
     \draw[black](M-1-4.north)--(M-6-4.south);
\end{tikzpicture}\raisebox{6.8ex}{.}
\hfill\hfill\hfill
\end{center}
The purposes of this paper is to prove the following Theorem.
\begin{theorem}[Sylvester's Determinant Identity]
For a square matrix $\AA$,
\begin{equation}\label{beforemain}
\det{\AA} \cdot \det{ \AA_{i,j|k,\ell}}
=
\det{ \AA_{i|k}} \cdot \det{\AA_{j|\ell }} - \det{\AA_{i|\ell}}\cdot \det{\AA_{j|k}}.
\end{equation}
\end{theorem}

\section{Proof}

Let us now prove \eqref{beforemain} by induction. As the sign of the determinant flips for each row or column that is permuted, it is clear that it is sufficient to show
\begin{equation}  \label{main}
\det{\AA}\cdot \det{ \AA_{1,2|1,2}} = 
\det{ \AA_{1|1}} \cdot \det{ \AA_{2|2}} -\det{ \AA_{1|2}} \cdot \det{ \AA_{2|1}}
\end{equation}
to prove  \eqref{beforemain} --- this simplifies the presentation somewhat.
\medskip

Recalling $|\AA| = 1$ when $\AA$ is $0 \times 0$, it is easy to verify \eqref{main} explicitly when $n=2$ 
$$
\det{
\begin{matrix}
a_{1,1} & a_{1,2} \\
a_{2,1} & a_{2,2}
\end{matrix}
}
\cdot
1
=
a_{2,2} \cdot a_{1,1} - a_{2,1} \cdot a_{1,2}
=
\det{ \AA_{1|1}} \cdot \det{ \AA_{2|2}} -\det{ \AA_{1|2}} \cdot \det{ \AA_{2|1}}.
$$

Extend $\AA$ by an extra row and column and denote this new $(n+1) \times (n+1)$ matrix by $\AA^{+}$. Further, let $\mathbb{A}^{(i)}$ be the matrix obtained by replacing the $i^{\rm th}$ row of $\AA$ with the last row of $\AA^{+}$ (the row that was added) less the `corner' element $a_{n+1,n+1}$.\medskip

We need to prove
\begin{equation}\label{induction}
\det{ \AA^{+}} \cdot \det{\AA_{1,2|1,2}^{+}} =\det{ \AA_{1|1}^{+}} \cdot \det{ \AA_{2|2}^{+}} -\det{ \AA_{1|2}^{+}} \cdot \det{ \AA_{2|1}^{+}}
\end{equation}
assuming the induction hypothesis (\ref{main}) holds.
\medskip


We have $\det{\AA}$ is a polynomial in $n^{2}$ variables 
with $n!$ terms comprised of a product of $n$ distinct $a_{i,j}$'s --- in other words $\det{\AA}$ is linear in each $a_{i,j}$. Therefore the $a_{i,j}$ cofactor can be written as  $\pdiff{~}{a_{i,j}}\det{\AA}$. Expanding the determinant along the (say) $j^{\rm th}$ column:\begin{equation}
\det{\AA} = \sum_{i=1}^{n}a_{i,j}\cdot (-1)^{i+j}\det{\AA_{i|j}} 
= \sum_{i=1}^{n}a_{i,j}\cdot \pdiff{\det\AA}{a_{i,j}}.
\end{equation}
Alternatively, using the $i^{\rm th}$ row gives
\begin{equation}
\det{\AA}=\sum_{j=1}^{n}a_{i,j}\cdot \frac{\partial\det{\AA}}{\partial a_{i,j}}
\end{equation}
and this implies that
\begin{equation}\label{Di}
\det{\AA^{(i)}} 
= \sum_{j=1}^{n}a_{n+1,j}\cdot \pdiff{\det{\AA}}{a_{i,j}}
\end{equation}

The last equivalence defines, for each $1 \leq i \leq n$, a linear differential operator $D^{(i)}$ with constant coefficients, namely
$$
D^{(i)}:= \sum_{j=1}^{n}a_{n+1,j}\cdot \pdiff{~}{a_{i,j}}.
$$
which obeys the product rule and commutes with $D^{(\ell )}$. Moreover
\begin{equation}
D^{(\ell )}\det{\AA^{(i)}} = 0
\end{equation}
because, when $\ell \neq i$, the LHS results in a determinant of a matrix having two identical rows (the $\ell^{\rm th}$ and $j^{\rm th}$). And, when $\ell=i$, the second $D^{(i)}$ is differentiating with respect to $a_{i,j}$ elements which are no longer part of the $D^{(i)}\det{\AA}$ polynomial.

\subsection{Main identity}

Expanding the determinant of $\AA^{+}$ along its last column, we get 
\begin{equation}\label{We8Jul1137}
\det{ \AA^{+}} =a_{n+1,n+1}\det{\AA}-\sum_{i=1}^{n}a_{i,n+1}D^{(i)}\det{\AA}
\end{equation}
since $-D^{(i)}\det{\AA}$ is now the cofactor of $a_{i,n+1}$.

It is important to realize that the last formula remains correct for \emph{each} determinant of (\ref{induction}) because (\ref{Di}) automatically `skips' (by contributing zero) the deleted elements.

\subsection{Proving (\ref{induction})}

Using \eqref{We8Jul1137} to expand the determinants of \eqref{induction} produces
\begin{align}
\MoveEqLeft\left( a_{n+1,n+1}\det{ \AA}-\sum_{i=1}^{n}a_{i,n+1}D^{(i)}\det{\AA}\right)\cdot \left( a_{n+1,n+1}\det{ \AA_{1,2|1,2}}-\sum_{i=1}^{n}a_{i,n+1}D^{(i)}\det{ \AA_{1,2|1,2}}\right)  \label{final} \\
=&~\left( a_{n+1,n+1}\det{ \AA_{1|1}}-\sum_{i=1}^{n}a_{i,n+1}D^{(i)}\det{ \AA_{1|1}}\right) \cdot \left( a_{n+1,n+1}\det{ \AA_{2|2}}-\sum_{i=1}^{n}a_{i,n+1}D^{(i)}\det{ \AA_{2|2}} \right) \notag\\
&-\left( a_{n+1,n+1}\det{ \AA_{1|2}}-\sum_{i=1}^{n}a_{i,n+1}D^{(i)}\det{ \AA_{1|2}}\right) \cdot \left( a_{n+1,n+1}\det{ \AA_{2|1}}-\sum_{i=1}^{n}a_{i,n+1}D^{(i)}\det{ \AA_{2|1}} \right) \notag
\end{align}

Collecting terms proportional to $a_{n+1,n+1}^{2}$ results in (\ref{main}); collecting terms
proportional to $a_{i,n+1}a_{j,n+1}$ results in
\begin{align*}
\MoveEqLeft \hspace{-2em}
D^{(i)}\det{\AA}\cdot D^{(j)}\det{ \AA_{1,2|1,2}} +D^{(j)}\det{\AA}\cdot D^{(i)}\det{ \AA_{1,2|1,2}} \\
=&~  D^{(i)}\det{ \AA_{1|1}} \cdot D^{(j)}\det{\AA_{2|2}} +D^{(j)}\det{ \AA_{1|1}} \cdot D^{(i)}\det{ \AA_{2|2}} \\
&\hspace{2em}-D^{(i)}\det{ \AA_{1|2}} \cdot D^{(j)}\det{ \AA_{2|1}} -D^{(j)}\det{ \AA_{1|2}}
\cdot D^{(i)}\det{ \AA_{2|1}} 
\end{align*}
which is the same as $D^{(i)}D^{(j)}$ applied to (\ref{main}), and thereby assumed correct (use the product rule twice and recall that $D^{(i)}D^{(j)}$ applied to a single determinant results in zero).


Finally, terms proportional to $a_{n+1,n+1}a_{i,n+1}$ yield
\begin{align*}
\MoveEqLeft \hspace{-2em} -D^{(i)}\det{\AA}\cdot \det{ \AA_{1,2|1,2}} -\det{\AA}\cdot
D^{(i)}\det{ \AA_{1,2|1,2}} \\
=& ~D^{(i)}\det{ \AA_{1|2}} \cdot \det{\AA_{2|1}} +\det{ \AA_{1|2}} \cdot D^{(i)}\det{ \AA_{2|1}} -D^{(i)}\det{ \AA_{1|1}} \cdot \det{ \AA_{2|2}} -\det{ \AA_{1|1}} \cdot D^{(i)}\det{ \AA_{2|2}}
\end{align*}
which is the same as $-D^{(i)}$ applied to (\ref{main}).\medskip

We have thus been able to cancel out all terms of \eqref{final}; the extended identity is thus verified.


\begin{thebibliography}{9}
\bibitem{akritas} A G Akritas, E K Akritas, G Malaschonok: ``Various poofs of
Sylvester's (determinant) identity'' \emph{Mathematics and Computers in
Simulation}, 1996, \textbf{42} \#4, 585-593
\end{thebibliography}
\end{document}